\documentclass[11pt,twoside]{amsart}

\title[Pimsner-Toeplitz algebras]{On the universal property of Pimsner-Toeplitz $C^*$-algebras and 
their continuous analogues}
\author{Ilan Hirshberg}
\address{Department of Mathematics, University of California, Berkeley, CA 
94720,USA}
\email{ilan@math.berkeley.edu}

\usepackage{amssymb}
\usepackage{amscd}
\usepackage{amsthm}

\theoremstyle{plain}

\newtheorem{Thm}{Theorem}[section]
\newtheorem{Cor}[Thm]{Corollary}
\newtheorem{Lemma}[Thm]{Lemma}

\newtheorem{Obs}[Thm]{Observation}

\theoremstyle{definition}
\newtheorem{Def}[Thm]{Definition}

\newtheorem{Rmk}[Thm]{Remark}

\newcommand{\B}{\mathcal{B}}
\newcommand{\A}{\mathcal{A}}

\newcommand{\W}{\mathcal{W}}

\newcommand{\T}{{\mathbb T}}
\newcommand{\R}{{\mathbb R}}
\newcommand{\N}{{\mathbb N}}
\newcommand{\Z}{{\mathbb Z}}
\newcommand{\C}{{\mathbb C}}
\newcommand{\EE}{\mathcal{E}}

\newcommand{\lb}{\left <}
\newcommand{\rb}{\right >}

\newcommand{\arrow}{\rightarrow}

\newcommand{\Ltwo}{L^2(\R_+)}

\newcommand{\LoneE}{L^1(E)}
\newcommand{\WE}{\W_E}
\newcommand{\TE}{\mathcal{T}_E}

\newcommand{\lnorm}{\left \|}
\newcommand{\rnorm}{\right \|}

\newcommand{\dint}{\int^{\oplus}}

\newcommand{\LtwoE}{\dint_{\R_+}E_xdx}

\newcommand{\psish}{\psi_{\#}}
\newcommand{\Hsh}{H_{\#}}
\newcommand{\eps}{\varepsilon}

\begin{document}

\begin{abstract}
We consider $C^*$-algebras generated by a single Hilbert bimodule 
(Pimsner-Toeplitz algebras) and by a product systems of Hilbert bimodules.
We give a new proof of a theorem of Pimsner, which states that any 
representation of the generating bimodule gives rise to a representation of 
the Pimsner-Toeplitz algebra. Our proof does not make use of the conditional 
expectation onto the subalgebra invariant under the dual action of the 
circle group.
We then prove the analogous statement for the case of product systems, 
generalizing a theorem of Arveson from the case of product systems of 
Hilbert spaces.
\end{abstract}
\maketitle

\section{Introduction}

Let $E$ be a Hilbert 
module over $\A$, equipped with a left action of $\A$ via adjointable 
operators. We shall refer to such an $E$ as a Hilbert bi-module. We assume that $E$ is full, i.e. $\overline{\lb E,E \rb} = \A$. We make no further assumptions on the left action of $\A$. 
Let $\B$ be a $C^*$-algebra. A 
\emph{covariant homomorphism} $\psi$ of $E$ into $\B$ is a $\C$-linear map 
$\psi_E : E \arrow \B$ along with a homomorphism $\psi_{\A}: \A \arrow \B$ 
such that for all $e,f\in E$, $a,b\in \A$ we have $\psi_E(aeb) = 
\psi_{\A}(a)\psi_E(e)\psi_{\A}(b)$, $\psi_E(e)^*\psi_E(f) = \psi_{\A}(\lb 
e,f\rb)$.
In the sequel, we will write $\psi$ for both $\psi_{\A}$ and $\psi_E$, when 
it causes no confusion.
When $\B = \B(H)$ for a Hilbert space $H$, we call a covariant homomorphism 
a \emph{representation}.

\begin{Rmk} In \cite{ms1} what we call a \emph{representation} is called an 
\emph{isometric covariant} representation. Since we won't deal with 
other kinds of representations considered in \cite{ms1}, we shall 
use the shorter terminology.
\end{Rmk}

We recall Pimsner's construction. Let $\EE = \bigoplus_{n=0}^{\infty} E^{\otimes n}$ (where we take 
$E^{\otimes 0} = \A$).
For $e \in E$, let $T_e \in \B(\EE)$ be given by $T_e(\xi) = e \otimes 
\xi$. The map sending $$e \mapsto T_e, \; \A \ni a \mapsto \textrm{ left 
multiplication by } a$$ is a covariant homomorphism of $E$ into $\B(\EE)$. 
We let $\TE$ be the $C^*$-subalgebra of $\B(\EE)$ generated by $\{T_e \; | 
\; e\in E\}$. We will lighten notation by identifying $e \in E$ with 
$T_e$, and $a$ with multiplication by $a$ in $\TE$, when it does not lead to 
confusion.

Note that if $\pi:\TE \arrow \B(H)$ is a representation, then the 
restrictions of $\pi$ to $E$ and $\A$ form a representation of $E$. Our goal 
in the first section will be is to give a new proof of the following theorem 
-- a restatement of a theorem of Pimsner (\cite{pimsner}, Theorem 3.4) -- 
which shows that any representation arises in this manner.
\begin{Thm}[Pimsner]
\label{pimsner universal thm} Let $\psi$ be a representation of $E$ on a 
Hilbert space $H$. The map $e \arrow \psi(e)$ extends to a homomorphism $\TE 
\arrow \B(H)$.
\end{Thm}

\begin{Rmk} Pimsner's proof relies on the conditional expectation map of the 
algebra $\TE$ onto the fixed point subalgebra for the dual action of $\T$, 
generalizing the proof for the Cuntz algebras from \cite{cuntz}. The 
motivation leading to the proof presented herein was to obtain the 
continuous analogue, Theorem \ref{cts univ property} below.
The continuous analogues, described below, admit an analogous action of 
$\R$, rather than $\T$. Thus, one cannot obtain a bounded expectation map by 
averaging the group action. We note that an unbounded expectation map has 
been used by Zacharias to study Arveson's spectral $C^*$-algebras in 
\cite{zacharias pi}.
We refer the reader to \cite{arveson} for more details on the spectral 
$C^*$-algebras, and to \cite{hz} for a recent survey.
\end{Rmk}

In \cite{hmps}, we considered a certain continuous analogue of 
the algebras $\TE$, generalizing to the context of Hilbert modules 
Arveson's spectral $C^*$-algebras (see \cite{arveson}). We recall the 
definitions.

\begin{Def} \label{meas. bundle} Let $\A$ be a separable $C^*$-algebra. A 
\emph{measurable bundle of Hilbert $\A$-bimodules over $\Omega$}, $E$, is a 
collection $\{E_x \; | \; x\in\Omega\}$ of right Hilbert $\A$-modules with 
left actions via adjointable operators, along with a distinguished vector 
subspace $\Gamma$ of $\Pi_{x\in\Omega}E_x$ (called the set of 
\emph{measurable sections}) such that
\begin{enumerate}
\item \label{meas. cond.} For any $\xi \in \Gamma$, $a\in\A$, the functions 
$x \mapsto \lb \xi(x),\xi(x)\rb$, $x \mapsto \lb a\xi(x),\xi(x)\rb$ are 
measurable (as functions $\Omega \mapsto \A$).
\item \label{closure cond.} If $\eta \in \Pi_{x\in\Omega}E_x$ satisfies that 
$x \mapsto \lb \xi(x),\eta(x)\rb$ is measurable for all $\xi\in\Gamma$ then 
$\eta\in \Gamma$.
\item \label{countable separation cond.} There exists a countable subset 
$\xi_1,\xi_2,...$ of $\Gamma$ such that for all $x\in \Omega$, 
$\xi_1(x),\xi_2(x),...$ are dense in $E_x$.
\end{enumerate}
\end{Def}
We refer the reader to the appendix of \cite{hmps} for more details.

\begin{Def} Let $\A$ be a separable $C^*$-algebra.
A \emph{product system of $\A$-bimodules} $E$ is a measurable bundle of 
$\A$-bimodules over $\R_+$, along with a multiplication map $E \times E 
\arrow E$, which descends to an isomorphism $E_x \otimes_{\A}E_y \arrow 
E_{x+y}$ for all $x,y\in \R_+$ (where $E_x$ is the fiber over $x$), and is 
measurable in the sense that if $\xi$ is a measurable section and $e \in 
E_y$ then the sections $x \mapsto e\xi(x-y)$, $x \mapsto \xi(x-y)e$ ($0$ if 
$x<y$) are also measurable.
\end{Def}

The elements $e\in E$ act on $\LtwoE$ on the left as adjointable operators, which we denote $W_e$, or by abuse of notation, just $e$.
Note that $\|e\|_{E_x} \geq \|W_e\|_{\B(\LtwoE)}$. 
Denote by 
$\LoneE$ the space of measurable sections $\xi$ that satisfy $\int_{\R_+} \|\xi (x)\| dx < \infty$.

\begin{Def}
For $f\in \LoneE$ we define $W_f \in \B\left (\LtwoE\right )$ by
$$
W_f = \int_{\R_+} W_{f(x)} dx
$$

We denote by $\W_E$ the $C^*$-subalgebra of $\B\left (\LtwoE\right )$ 
generated by
$$\{W_f \; | \; f\in\LoneE\}$$
\end{Def}
We refer the reader to \cite{hmps} for examples, and a discussion of the 
$K$-theory of $\WE$.

\begin{Def} Let $E$ be a product system of $\A$-bimodules. A 
\emph{representation} $\psi$ of $E$ on $H$ is a map $\psi_E : E \arrow 
\B(H)$, along with a representation $\psi_{\A}: \A \arrow \B(H)$ such that
\begin{enumerate}
\item The restriction of $\psi$ to each fiber of $E$ is a representation of 
the fiber.
\item For any $e,f \in E$, $\psi(e)\psi(f) = \psi(ef)$.
\item If $\xi$ is a measurable section of $E$ then $x \mapsto \psi(\xi(x))$ 
is a weakly measurable function.
\item $\overline{\bigcup_{x>0}\psi(E_x)H} = \overline{\psi(A)H}$.
\end{enumerate}
\end{Def}

If $x \mapsto f(x)$ is a measurable section of $E$ satisfying 
$\int_{\R_+} \|f(x)\|dx < \infty$ and $\psi$ is a representation of $E$ on $H$, then 
we have an integrated form of the representation
$$\psi(f) = \int_{\R_+}\psi(f(x))dx$$

Our goal in the second part of this paper will be to prove the following 
continuous analogue of Theorem \ref{pimsner universal thm}.

\begin{Thm}
\label{cts univ property}
Let $\psi$ be a representation of $E$ on a Hilbert space $H$. The map $W_f 
\arrow \psi(f)$ extends to a homomorphism $\WE \arrow \B(H)$.
\end{Thm}
This theorem generalizes a theorem of Arveson (\cite{arveson}, Theorem 
4.6.6) from the case of product systems of Hilbert spaces. Specializing our proof below to the case of Hilbert spaces will give a simpler approach to Arveson's theorem.

\section{The discrete case -- proof of Theorem \ref{pimsner universal thm}}

\begin{Def} 
\label{Def: majorization}
Let $\psi,\rho$ be two representations of $E$ on $H$. We say 
that $\psi$ \emph{majorizes} $\rho$, and write $\psi \succ \rho$, if for any 
$e_1,...,e_n \in E$ and any polynomial $p$ in $2n$ non-commuting variables, 
we have
$$
\lnorm p(\psi(e_1),...,\psi(e_n),\psi(e_1)^*,...,\psi(e_n)^*) \rnorm \geq 
\lnorm p(\rho(e_1),...,\rho(e_n),\rho(e_1)^*,...,\rho(e_n)^*) \rnorm
$$
In other words, $\psi \succ \rho$ if there is a (necessarily unique) 
homomorphism $C^*(\{\psi(e) \; | \; e \in E\}) \arrow C^*(\{\rho(e) \; | \; 
e \in E\})$ which satisfies $\psi(e) \mapsto \rho(e)$ for all $e\in E$.

If $\psi \succ \rho$ and $\psi \prec \rho$, we write $\psi \approx \rho$.

We say that $T \succ \psi$ if the map $T_e \arrow \psi(e)$ extends to a 
homomorphism $\TE \arrow C^*(\{\psi(e) \; | \; e \in E\})$, i.e. if
$$
\lnorm p(\psi(e_1),...,\psi(e_n),\psi(e_1)^*,...,\psi(e_n)^*) \rnorm \leq 
\lnorm p(e_1,...,e_n,e_1^*,...,e_n^*) \rnorm_{\B(\EE)}
$$
\end{Def}
Thus Theorem \ref{pimsner universal thm} states that $T \succ \psi$ for any 
representation $\psi$ of $E$.

Note that the relation $\succ$ is clearly transitive.

We first recall the following lemma (noted in \cite{ms1} and in references 
therein). The proof is straightforward.
\begin{Lemma} \label{Lemma: contraction is isometry} Let $\psi$ be a 
representation of $E$ on $H$. Regarding $H$ as a right $\A$-module via 
$\psi$, we form the tensor product $E \otimes_{\A} H$ to obtain a Hilbert 
space. The contraction map $$e \otimes \xi \mapsto \psi(e)\xi \;\;\;\; e \in 
E\,, \; \xi \in H$$ is well defined and extends to an isometry $$E 
\otimes_{\A} H \arrow H$$
\end{Lemma}

If $\psi$ is a representation of $E$, $n >0$, then we can define a 
representation of $E^{\otimes n}$ by $e_1 \otimes \cdots \otimes e_n \mapsto 
\psi(e_1)\psi(e_2) \cdots \psi(e_n)$. We will denote this representation by 
$\psi$ as well.

We may assume without loss of generality that $\psi_{\A}$ is non-degenerate, 
and we make this assumption throughout (i.e., we assume throughout 
that $\overline{\psi(\A)H} = H$).

If $\pi$ is a representation of $\A$ on a Hilbert space $H$, then we can 
define a representation $T\otimes_{\A}1$ of $E$ on $\EE \otimes_{\pi}H$ 
(this is called an \emph{induced} representation in \cite{ms1}). We note 
that clearly, $T \succ T\otimes_{\A}1$.
The following lemma generalizes the fact that any isometry $S$ which 
satisfies $S^nS^{n*} \arrow 0$ in the strong operator topology is unitarily 
equivalent to a direct sum of copies of the unilateral shift on $\ell^2$. The 
reader can find a proof in \cite{ms1}.

\begin{Lemma}
\label{Wold lemma}
Let $\psi$ be a representation of $E$ on $H$, such that 
$\bigcap_{n>0}\overline{\psi(E^{\otimes n})H} = \{0\}$. let $H_0 = 
(\psi(E)H)^{\perp}$.
\begin{enumerate}
\item $H_0$ is invariant for $\psi(\A)$.
\item Let $H_n = \overline{\psi(E^{\otimes n})H_0}$. We have $H = 
\bigoplus_{n=0}^{\infty}H_n$.
\item For any $n$, $W_n : E^{\otimes n} \otimes_{\A} H_0 \arrow H_n$ given 
by the contraction
$$
W_n(e_1 \otimes \cdots \otimes e_n \otimes \xi) = \psi(e_1)\psi(e_2)\cdots 
\psi(e_n)\xi
$$
is a well defined unitary operator (where for $W_0$ is the contraction $a 
\otimes \xi \mapsto \psi(a)\xi$, $a\in E^{\otimes 0} = \A$).
\item $W = \bigoplus_{n = 0}^{\infty}W_n : \EE \otimes_{\A} H_0 \arrow H$ is 
a unitary operator which satisfies
$$
W(T_e \otimes_{\A} 1) = \psi(e)W
$$
for all $e \in E$, i.e. it implements a unitary equivalence between the 
covariant representations $T \otimes_{\A} 1_{H_0}$ and $\psi$.
\end{enumerate}
\end{Lemma}

\begin{Cor}
\label{singular reduction}
Let $\psi$ be as in Lemma \ref{Wold lemma}, then $T \succ \psi$.
\end{Cor}

Now let $\psi$ be any (non-degenerate) representation of $E$ on $H$. By 
Corollary \ref{singular reduction}, to prove Theorem \ref{pimsner universal 
thm}, it suffices to show that $\psi$ is majorized by a representation which satisfies the condition of Lemma \ref{Wold lemma}.

For any $\lambda \in \T$, we define a representation $\psi_{\lambda}$, given 
by $\psi_{\lambda}(e) = \lambda\psi(e)$, $\psi_{\lambda}(a) = \psi(a)$, 
$e\in E$, $a\in \A$.
We can now form a direct integral to obtain a representation 
$\bar{\psi}$ on $H \otimes L^2(\T)$, given by
$$
\bar{\psi} = \int_{\T}^{\oplus}\psi_{\lambda}d\lambda
$$
Since $\psi_{\lambda}(e) \arrow \psi(e)$ as $\lambda \arrow 1$ for all $e$ (in norm), we can 
easily see that $\bar{\psi} \succ \psi$.

Let $U$ be the bilateral shift on $\ell^2(\Z)$. We form a representation 
$\tilde{\psi} = \psi \otimes U$ of $E$ on $H \otimes \ell^2(\Z)$ by $\tilde{\psi}(e) = \psi(e) \otimes U$, $\tilde{\psi}(a) = a \otimes 1$, 
$e\in E$, $a\in \A$.
Applying the Fourier transform to the second variable shows that
$\bar{\psi}$ and $\tilde{\psi}$ are unitarily equivalent. (i.e. they are 
intertwined by a unitary). In particular, we have $\bar{\psi} \approx 
\tilde{\psi}$.

Denote by $P_+$ the projection of $H \otimes \ell^2(\Z)$ onto $H \otimes 
\ell^2(\N)$.
We denote by $\psi_+$ the restriction of $\tilde{\psi}$ to the invariant 
subspace $H \otimes \ell^2(\N)$, i.e. $\psi_+(e) = \tilde{\psi}(e)P_+$ (where here we will think of $\psi_+$ as both a representation on $H \otimes \ell^2(\N)$ and as a (degenerate) representation on $H \otimes \ell^2(\Z)$).
 Denote 
by $S$ the unilateral shift on $\ell^2(\N)$, and let $V = 1_H \otimes S$.

\begin{Obs} \label{4.5.3} For any $k$, any polynomial $p(x_1,...,x_{2k})$ 
in $2k$ non-commuting variables and any $e_1,...,e_k\in E$,
and any $m>\deg (p)$, we have
$$
V^{m*} p(\psi_+(e_1),...,\psi_+(e_k)^*)  V^m - P_+ 
p(\tilde{\psi}(e_1),...,\tilde{\psi}(e_k)^*) P_+ = 0
$$
\end{Obs}
We leave the straightforward verification to the reader.

\begin{Lemma} \label{4.5.4} If $A$ is in the $*$-algebra generated by 
$$\{\tilde{\psi}(e) \, | \, e\in E\}$$ then $$\|P_+AP_+\| = \|A\|$$
\end{Lemma}
\begin{proof}
Note that any operator of the form $\tilde{\psi}(e)$ commutes with all 
operators of the form $1_H \otimes U^n$, $n \in \Z$. Therefore $A$ commutes 
with $1\otimes U^n$, $n\in\Z$ as well. Let $P_n$ denote the projection onto 
$H \otimes \ell^2(\{n,n+1,...\})$ (so $P_0 = P_+$). We have $(1 \otimes U^m) 
P_n (1 \otimes U^m)^* = P_{n+m}$ for all $n,m\in\Z$, and therefore we have
$$
(1 \otimes U^n) P_+ A P_+ (1 \otimes U^n)^* = P_{n} A P_{n}
$$
so $\|P_+AP_+\| = \|P_n A P_n\|$ for all $n \in \Z$. Since $P_n \arrow 
1_{\ell^2(\Z)}$ as $n \arrow - \infty$ in the strong operator topology, we have
$$ \|P_+AP_+\| = \lim_{n \arrow -\infty} \|P_n A P_n\| = \|A\| $$
as required.
\end{proof}

\begin{Cor} $\psi_+ \succ \tilde{\psi}$.
\end{Cor}
\begin{proof}
Let $e_1,...,e_k \in E$, and let $p(x_1,...,x_{2k})$ be a polynomial in $2k$ 
non-commuting variables. Since the $V_m$ are isometries, we have
$$
\left \|V^{m*}p(\psi_+(e_1),...,\psi_+(e_k)^*) V^m \right \|
\leq
\left \|p(\psi_+(e_1),...,\psi_+(e_k)^*) \right \|
$$
so by Observation \ref{4.5.3}, we have $$
\left \| P_+ p(\tilde{\psi}(e_1),...,\tilde{\psi}(e_k)^*) P_+  \right \|
\leq
\left \|p(\psi_+(e_1),...,\psi_+(e_k)^*) \right \|
$$
and by Lemma \ref{4.5.4},
$$
\left \| P_+ p(\tilde{\psi}(e_1),...,\tilde{\psi}(e_k)^*) P_+  \right \|
=
\left \|  p(\tilde{\psi}(e_1),...,\tilde{\psi}(e_k)^*)  \right \|
$$
\end{proof}

\begin{proof}[Proof of Theorem \ref{pimsner universal thm}] Note that 
$\psi_+$ satisfies the conditions of Lemma \ref{Wold lemma}. It 
therefore suffices to show that $\psi_+ \succ \psi$, and indeed, we saw that 
$\psi_+ \succ \tilde{\psi}$ and $\tilde{\psi} \succ \psi$.
\end{proof}

\begin{Rmk} The proof in this section was obtained in the course of the 
author's dissertation work under the supervision of W.B. Arveson, and is 
motivated by ideas from \cite{arveson}.
\end{Rmk}

\section{The continuous case -- proof of Theorem \ref{cts univ property}}
The approach here will differ from the proof above for the discrete case, in that we do not have a continuous analogue of Lemma \ref{Wold lemma} (see Remark \ref{Cooper remark} below). Aside for that, we shall follow a similar path.

We begin by giving the analogue of Definition \ref{Def: majorization}.
\begin{Def} 
Let $\psi,\rho$ be two representations of a product system $E$ (over $\A$) on $H$. We say 
that $\psi$ \emph{majorizes} $\rho$, and write $\psi \succ \rho$, if for any 
$f_1,...,f_n \in L^1(E)$ and any polynomial $p$ in $2n$ non-commuting variables, 
we have
$$
\lnorm p(\psi(f_1),...,\psi(f_n),\psi(f_1)^*,...,\psi(f_n)^*) \rnorm \geq 
\lnorm p(\rho(f_1),...,\rho(f_n),\rho(f_1)^*,...,\rho(f_n)^*) \rnorm
$$
In other words, $\psi \succ \rho$ if there is a (necessarily unique) 
homomorphism $C^*(\{\psi(f) \; | \; f \in L^1(E)\}) \arrow C^*(\{\rho(f) \; | \; 
f \in L^1(E)\})$ which satisfies $\psi(f) \mapsto \rho(f)$ for all $f\in L^1(E)$.

If $\psi \succ \rho$ and $\psi \prec \rho$, we write $\psi \approx \rho$.

We say that $W \succ \psi$ if the map $W_f \arrow \psi(f)$ extends to a 
homomorphism $\WE \arrow C^*(\{\psi(f) \; | \; f \in L^1(E)\})$.
\end{Def}
As in the discrete case, the relation $\succ$ is clearly transitive. Theorem \ref{cts univ property} states that $W \succ \psi$ for any 
representation $\psi$ of $E$.

As in the discrete case, we may assume without loss of generality that $\psi_{\A}$ is non-degenerate, 
and we make this assumption throughout.

\begin{Def} Let $\psi$ be a representation of $E$ on $H$. A subspace $H'$ of $H$ is said to be \emph{invariant} for $\psi$ if $\psi(E_x)H' \subseteq H'$ (for all $x>0$) and $\psi(\A)H' \subseteq H'$. $H'$ will be said to be \emph{reducing} if it is invariant, and furthermore $\psi(e)^*H' \subseteq H'$ for all $e \in E$.
\end{Def}
Let $\psi$ be a representation of $E$ on $H$, and let $H'$ be invariant for $\psi$, then we have a 
representation of $E$ on $H'$ by restriction. Let $P$ be the projection onto $H'$, and let $\psi'$ denote the restriction, then $\psi'(f) = \psi(f)P$. Notice that if $H'$ is furthermore reducing, then $\psi \succ \psi'$.

We will make use of the following approximation lemma. The proof is straightforward, and we leave it to the reader.

\begin{Lemma} \label{approximation lemma}
 Let $\psi$ be a representation of $E$ on $H$. Suppose that there is a sequence of projections $P_n \arrow 1$ in the strong operator toplogy, such that $P_nH$ is invariant for $\psi$ for all $n$. Denote by $\psi_n$ the restricted representation of $\psi$ to $P_nH$. For any polynomial $p(x_1,...,x_{2k})$ in $2k$ non-commuting variables and $f_1,...,f_k\in \LoneE$, if $\lnorm p(\psi_n(f_1),...,\psi_n(f_k),\psi_n(f_1)^*,...,\psi_n(f_k)^*) \rnorm \leq M$ for all $n$ (for some constant $M$), then $\lnorm p(\psi(f_1),...,\psi(f_k),\psi(f_1)^*,...,\psi(f_k)^*) \rnorm \leq M$. 

Consequently, if $\rho$ is a representation of $E$ such that $\rho \succ \psi_n$ for all $n$ then $\rho \succ \psi$.
\end{Lemma}

If $\pi$ is a representation of $\A$ on $H$, we can form a representation $W\otimes_{\A}1$ of $\WE$ on $\LtwoE \otimes_{\pi}H$, as in the discrete case. We clearly have $W \succ W \otimes_{\A}1$. 

Let $S_x:L^2(\R_+) \arrow L^2(\R_+)$ denote the unilateral shift semigroup. We form a representation $\psi_+$ on $H \otimes L^2(\R_+)$ by $\psi_+(e) = \psi(e) \otimes S_x$ ($e \in E_x$), $\psi_+(a) = \psi(a) \otimes 1$. The following Lemma is a straightforward generalization, with a small improvement, of an argument in the proof of \cite{arveson}, Theorem 4.4.3. We include a full proof for the reader's convenience. 
\begin{Lemma}
\label{cts 4.4.3}
\

\begin{enumerate}
\item Let $\psi$ be a (non-degenerate) representation of $E$ on a Hilbert space $H$. There is a unique isometry 
$$C: \LtwoE \otimes_{\A}H \arrow H\otimes\Ltwo \cong L^2(\R_+,H)$$ 
satisfying 
$$
C(f \otimes \xi)(x) = (\psi(f(x))\xi)
$$
for any $\xi \in H$, $f \in \LoneE$ such that $\int_{\R_+}\|f(x)\|^2dx < \infty$ (those are dense both in $\LoneE$ and in $\LtwoE$).

The range of $C$ is $H_{\#} = \{\xi\in L^2(\R_+,H) \; | \; \xi(x)\in \overline{\psi(E_x)H} \textrm{ a.e.}\, x\}$.

\item 
$H_{\#}$ is invariant for the representation $\psi_+$. Denote the restriction of $\psi_+$ to $H_{\#}$ by $\psi_{\#}$. $C$ implements a unitary equivalence between $\psi_{\#}$ and $W \otimes_{\A} 1$ of $E$ on $\LtwoE \otimes_{\A} H$. 
\end{enumerate}
\end{Lemma}
\begin{proof}
$C$ is well defined on the given domain, which is total in $\LtwoE \otimes_{\A}H$. To show that $C$ extends to an isometry, it suffices to check inner products on those vectors. So, for $f,g \in \LoneE$ such that $\int_{\R_+}\|f(x)\|^2dx<\infty$,$\int_{\R_+}\|g(x)\|^2dx < \infty$, and $\xi,\eta\in H$, we have:
$$
\lb \psi(f(x))\xi,\psi(g(x))\eta \rb_{H} 
= 
\lb \xi,\psi(f(x))^*\psi(g(x))\eta \rb_{H} 
=
\lb \xi,\psi(\lb f(x),g(x)\rb_{\A})\eta \rb_{H} = 
$$
$$
= \lb \psi(f(x))\otimes_{\A}\xi,\psi(g(x))\otimes_{\A}\eta \rb_{E_x \otimes_{\A} H}
$$
and now integrating both sides $dx$ gives the required identity. 

To prove that $C$ has the required range, we first note that the range of $C$ is clearly contained in $H_{\#}$. For the converse, suppose $g\in H_{\#}$ is orthogonal to the range of $C$. We must prove that $g = 0$. 

Let $e_1(x),e_2(x),...$ be a sequence of measurable sections of $E$ such that for all $x$, $\{e_1(x),e_2(x),...\}$ are total in $E_x$ (we are guaranteed the existence of such sequence by the definition of a product system). We may assume that all those sections are bounded, without loss of generality.
Let $u_1,u_2,...$ be a dense sequence in $L^1(\R_+) \cap L^2(\R_+)$. Let $\xi_1,\xi_2,...$ be a dense sequence in $H$. So, for all $m,n,p$, we have:
$$
\int_{0}^{\infty} u_m(x)\lb \psi(e_n(x))\xi_p,g(x)\rb = 0
$$
Since $u_1,u_2,...$ are dense in $L^2(\R_+)$, we see that $\lb \psi(e_n(x))\xi_p,g(x)\rb = 0$ a.e. $x$. Therefore, we have for all $n,p$ and a.e. $x$, $\lb \psi(e_n(x))\xi_p,g(x)\rb = 0$, and therefore, $g(x) = 0$ a.e., as required.

Finally, we need to check that $C(W_e \otimes_{\A} 1) = (\psi(e) \otimes S_x) C$ for all $e \in E_x$, $x \in \R_+$. It suffices to check this for vectors of the form $f \otimes_{\A} \xi$ as in the statement. Indeed, $C(W_e \otimes_{\A} 1) (f \otimes_{\A} \xi)(y) = C (ef\otimes_{\A}\xi)(y) = \psi((e\cdot f)(y))\xi = \psi(e\cdot f(y-x))\xi = \psi(e)\psi(f(y-x))\xi = (\psi(e)\otimes S_x) (C(f\otimes_{\A}\xi))$, as required (where $f(y-x)$ is understood to mean $0$ if $x>y$).
\end{proof}

\begin{Lemma} 
\label{quadrant-octant lemma}
Let $\psi$ be a representation of $E$ on $H$. Let $\psi_+$, $\psi_{\#}$, $H_{\#}$ be as in Lemma \ref{cts 4.4.3}, then $\psi_+ \approx \psish$.
\end{Lemma}
\begin{proof}
Since $\Hsh$ is a reducing subspace, we have $\psi_+ \succ \psish$.
Thus it remains to show that $\psish \succ \psi_+$.
By Lemma \ref{approximation lemma}, it suffices to exhibit projections $P_{\varepsilon} \in \B(L^2(\R_+,H))$ such that $P_{\varepsilon}(L^2(\R_+,H)$ is invariant for $\psi_+$, $P_{\eps} \arrow 1$ as $\eps \arrow 0$ (in the strong operator topology), and $\psish \succ \psi_{\eps}$ where $\psi_{\eps}$ denotes the restriction of $\psi_+$ to $P_{\varepsilon}(L^2(\R_+,H)$.

Denote $H_x = \overline{\psi(E_x)H}$.
Let $K_{\eps} = \{\xi \in L^2(\R_+,H) \; | \; \xi(x)\in H_x \ominus H_{x+\eps}\} \subseteq \Hsh$. $K_{\eps}$ is reducing for $\psish$ (and for $\psi_+$). Denote the restriction of $\psish$ to $K_{\eps}$ by $\psish^{\eps}$.

Now, for $n=1,2,...$, let $K_{\eps}^n =  \{\xi \in L^2(\R_+,H) \; | \; \xi(x)\in H_{x-n\eps} \ominus H_{x-(n-1)\eps}\}$. Note that the $K_{\eps}^n$ are mutually orthogonal, and are all orthogonal to $\Hsh$.

$K_{\eps}^n$ is invariant (but not reducing) for $\psi_+$. Note that if $\xi \in K_{\eps}^n$ then $\xi(x) = 0$ for a.e. $x\leq n\eps$. Let $\psi_{\eps}^n$ denote the restriction of $\psi_+$ to $K_{\eps}^n$. 

Define $U^n_{\eps}:K_{\eps} \arrow K_{\eps}^n$ by $U^n_{\eps}(\xi)(x) = \xi(x-n\eps)$ (where $\xi(x-n\eps)$ is understood to be $0$ if $x\leq n\eps$). It is easy to check that $U^n_{\eps}$ is unitary, and implements a unitary equivalence between $\psish^{\eps}$ and $\psi_{\eps}^n$. 

Let $H_+^{\eps} = \Hsh \oplus \bigoplus_{n=1}^{\infty}K_{\eps}^n \subseteq L^2(\R_+,H)$. Note that this space is invariant for $\psi_+$. Let $P_{\eps}$ be the projection onto the $H_+^{\eps}$, and $\psi_{\eps}$ the restriction of $\psi_+$ to $H_+^{\eps}$. So $P_{\eps} \arrow 1$ (since, for example, the range of $P_{\eps}$ contains $L^2(\R_+,H_{\eps})$), and $\psish \succ \psi_{\eps}$ for all $\eps$, which is what we needed.  
\end{proof}

We may now proceed as in the discrete case.
For $\lambda \in \R$, we define a representation $\psi_{\lambda}$  
by $\psi_{\lambda}(e) = e^{ix\lambda}\psi(e)$, $\psi_{\lambda}(a) = \psi(a)$. We form a representation 
$\bar{\psi}$ on $H \otimes L^2(\R)$, by
$$
\bar{\psi} = \int_{\R}^{\oplus}\psi_{\lambda}d\lambda
$$
and since $\psi_{\lambda}(e) \arrow \psi(e)$ as $\lambda \arrow 0$ for all $e$ (in norm), we have $\bar{\psi} \succ \psi$.

Let $U_x$ be the bilateral shift group on $L^2(\R)$. We form a representation 
$\tilde{\psi}$ of $E$ on $H \otimes L^2(\R)$ by $\tilde{\psi}(e) = \psi(e) \otimes U_x$ ($e\in E_x$), $\tilde{\psi}(a) = a \otimes 1$.
Using the Fourier transform, we see that $\bar{\psi}$ and $\tilde{\psi}$ are unitarily equivalent, so $\tilde{\psi} \succ \psi$.

Let $V_x = 1_H \otimes S_x$, and let $P_+$ the projection of $H \otimes L^2(\R)$ onto $H \otimes \Ltwo$

The following are immediate analogues of Observation \ref{4.5.3} and Lemma \ref{4.5.4} above (and immediate generalizations of 4.5.3 and 4.5.4 in \cite{arveson}). We leave the simple proofs to the reader.

\begin{Obs}\label{cts 4.5.3}
For any $k$ and any polynomial $p(x_1,...,x_{2k})$ 
in $2k$ non-commuting variables and any $f_1,...,f_k\in L^1(E)$,
We have
$$
\lim_{x \arrow \infty} \lnorm V_x^* p(\psi_+(f_1),...,\psi_+(f_k)^*)  V_x - P_+p(\tilde{\psi}(f_1),...,\tilde{\psi}(f_k)^*) P_+ \rnorm = 0
$$
\end{Obs}

\begin{Lemma}\label{cts 4.5.4} If $A$ is in the $*$-algebra generated by 
$$\{\tilde{\psi}(f) \, | \, f\in L^1(E)\}$$ then $$\|P_+AP_+\| = \|A\|$$
\end{Lemma}

\begin{Cor} \label{4.5.3-4} $\psi_+ \succ \tilde{\psi}$.
\end{Cor}

\begin{proof}[Proof of Theorem \ref{cts univ property}]
Let $\psi$ be a (non-degenerate) representation of $E$ on $H$. We want to show that $W \succ \psi$. So, $W \succ W \otimes_{\A} 1$, $W \otimes_{\A} 1 \approx \psish$ (by Lemma \ref{cts 4.4.3}),  $\psish \approx \psi_+$ (by Lemma \ref{quadrant-octant lemma}), $\psi_+ \succ \tilde{\psi}$ (Corollary \ref{4.5.3-4}), and $\tilde{\psi} \succ \psi$ (as remarked above), concluding the argument.
\end{proof}

\begin{Rmk}
\label{Cooper remark}
There is a continuous analogue of the Wold decomposition, due to Cooper (\cite{cooper}), which states that if $S_x$ ($x>0$) is a strongly continuous semigroup of isometries on a Hilbert space, then $S_x$ is unitarily equivalent to a direct sum of a one-parameter unitary group and copies of the unilateral shift semigroup on $L^2(\R_+)$. 

Unlike the case of a single bimodule, this does not quite generalize to product systems. There is an approximate version, due to Arveson (\cite{arveson paper 3, arveson} for product systems of Hilbert spaces. Arveson's theorem states the following. Let $E$ is a product system of Hilbert \emph{spaces}, and let $\psi$ be a representation of $E$ on $H$ such that $\bigcap_{x>0}\overline{\psi(E_x)H} = \{0\}$. For any $\varepsilon>0$, let $H_{\varepsilon} = \overline{\psi(E_{\varepsilon})H}$, and let $\psi_{\varepsilon}$ be the restriction of $\psi$ to $H_{\varepsilon}$, then $\psi_{\varepsilon}$ is unitarily equivalent to a direct sum of the regular representation of $E$ (on $\LtwoE$, which is a Hilbert space here).
However, Arveson showed in \cite{arveson paper 3} that $H_{\varepsilon}$ cannot be replaced by $H$ in the theorem. This theorem was used by Arveson to prove the special case of Theorem \ref{cts univ property} for Hilbert spaces.
We do not know if the generalization of Arveson's theorem to the case of Hilbert modules holds.
\end{Rmk}

\end{document}